\documentclass{article}
\usepackage[utf8]{inputenc}
\usepackage{amsmath, amsfonts,amsthm}
\usepackage[left=1in,top=1in,right=1in,bottom=1in,letterpaper]{geometry}
\usepackage{multirow,booktabs,makecell}
\usepackage{graphicx,epstopdf,subcaption}
\usepackage{tikz}
\usetikzlibrary{positioning}
\usepackage{color,hyperref}
\usepackage{algorithm,algorithmic}
\usepackage[numbers]{natbib}

\newcommand{\bbR}{\mathbb{R}}
\newcommand{\calP}{\mathcal{P}}

\newcommand{\bbT}{\mathbb{T}}

\newcommand{\dxp}{D_x^+}
\newcommand{\dxm}{D_x^-}
\newcommand{\dxc}{D_x^c}

\newcommand{\myij}{{i,n}}
\newcommand{\ipj}{{i+1,n}}
\newcommand{\imj}{{i-1,n}}
\newcommand{\ijp}{{i,n+1}}

\newcommand{\dx}{\Delta x}
\newcommand{\dt}{\Delta t}

\newcommand{\hamlf}{H^{\text{LF}}}
\newcommand{\pxp}{p^+}
\newcommand{\pxm}{p^-}

\newcommand{\calG}{\mathcal{G}}

\newcommand{\dd}{\mathrm{d}}

\DeclareMathOperator*{\argmin}{\arg\min}

\newtheorem{theorem}{Theorem}[section]

\newtheorem{proposition}[theorem]{Proposition}

\theoremstyle{definition}

\newtheorem{problem}[theorem]{Problem}

\theoremstyle{remark}
\newtheorem{remark}[theorem]{Remark}

\title{
Equilibrium Correction Iteration for \\A Class of Mean-Field Game Inverse Problems
}
\author{Jiajia Yu \and Jian-Guo Liu \and Hongkai Zhao}
\date{}

\begin{document}

\maketitle

\begin{abstract}

This work investigates the ambient potential identification problem in inverse Mean-Field Games (MFGs), where the goal is to recover the unknown potential from the value function at equilibrium.
We propose a simple yet effective iterative strategy, Equilibrium Correction Iteration (ECI), that leverages the structure of MFGs rather than relying on generic optimization formulations. ECI uncovers hidden information from equilibrium measurements, offering a new perspective on inverse MFGs. 
To improve computational efficiency, two acceleration variants are introduced: Best Response Iteration (BRI), which uses inexact forward solvers, and Hierarchical ECI (HECI), which incorporates multilevel grids. While BRI performs efficiently in general settings, HECI proves particularly effective in recovering low-frequency potentials.
We also highlight a connection between the potential identification problem in inverse MFGs and inverse linear parabolic equations, suggesting promising directions for future theoretical analysis.
Finally, comprehensive numerical experiments demonstrate how viscosity, terminal time, and interaction costs can influence the well-posedness of the inverse problem.

\end{abstract}

\section{Introduction}
\label{sec: intro}

A Mean-Field Game (MFG) \cite{lasry2007mean,huang2006large,huang2007large} is a non-cooperative game involving a continuum of indistinguishable players. 
The individual player seeks the best response to the population state distribution, and meanwhile, the aggregate of individual best responses determines the population state distribution.
Solving an MFG entails finding a Nash equilibrium for the system, where no individual player can unilaterally improve their strategy. 
In other words, if the population adopts the Nash equilibrium strategy, the best response for any individual player is to also follow the Nash equilibrium strategy.
In general, an MFG system is a coupled PDE system consisting of a time-backward Hamilton-Jacobi-Bellman (HJB) equation for the value function $\phi$ and a time-forward Fokker-Planck (FP) equation for the distribution $\rho$ as follows:
\begin{equation}\left\{
\begin{aligned}
    -\partial_t\phi(x,t)-\nu\Delta\phi(x,t)+\tilde{H}(x,\nabla\phi(x,t),\rho(\cdot,t))=0,
    \quad &(x,t)\in \bbT^d\times[0,T],  \\
    \partial_t\rho(x,t)-\nu\Delta\rho(x,t)-\nabla\cdot(\rho D_p\tilde{H}(x,\nabla\phi,\rho))(x,t)=0,
    \quad &(x,t)\in \bbT^d\times[0,T],  \\
    \phi(x,T)=f_T(x,\rho(\cdot,T)),\quad \rho(x,0)=\rho_0(x), 
    \quad &x\in\bbT^d.    
\end{aligned}\right.
\label{eq: mfg general}
\end{equation}
Here $\tilde{H}$ is the Hamiltonian defined as the conjugate of a cost $\tilde{L}(x,v,\rho)$,
\begin{equation*}
    \tilde{H}:\Omega\times\bbR^d\times\calP(\bbT^d)\to\bbR, 
    \quad 
    \tilde{H}(x,p,\rho):=\sup_{v} ~ -\langle p,v \rangle - \tilde{L}(x,v,\rho).
\end{equation*}
The dependence of $\tilde{L}$ on $v$ reflects the dynamic cost, and the dependence on $\rho$ reflects the interaction cost with the population.
A solution $(\rho,\phi)$ to \eqref{eq: mfg general} is called a Mean-Field Nash Equilibrium (MFNE).
The forward problem is to solve for the Nash equilibrium $(\rho,\phi)$ given the cost function $\tilde{L}$.

While MFG has been widely applied in various fields, in many scenarios, the cost function is not directly observable, and only the Nash equilibrium can be measured. For instance, in the stock market, the Nash equilibrium corresponds to the equilibrium price which is observable, but the cost function represents the utility function of the traders which is not clearly known. 

This motivates the study of the inverse problem, which aims to identify the cost function $\tilde{L}$ based on measurements of the Nash equilibrium. 
In this paper, we focus on a separable Hamiltonian of the form
\begin{equation*}
    \tilde{H}(x,p,\rho) = H(p) - q(x) - f(x,\rho).
\end{equation*}
The inverse problem we address is to identify $q$ from the measurement of the Nash equilibrium $\phi(\cdot,0)$.
The separable Hamiltonian is induced from a separable cost function of the form
\begin{equation*}
    \tilde{L}(x,v,\rho) = L(v) + q(x) + f(x,\rho),
\end{equation*}
where $L$ represents the dynamic cost, $f$ is the interaction cost, and $q$ is referred to as the ambient potential.
The ambient potential $q(x)$ quantifies the cost associated with an individual state $x$ and is determined by the environment, independent of the population. 
In much of the inverse MFG literature, $q$ is often referred to as the obstacle. However, this terminology can be misleading, as we consider a general form of $q$, not necessarily a discontinuous or jump function.
A more appropriate term is potential, as used in many physical models. To avoid confusion with the notion of potential in potential mean-field games, we adopt the term ambient potential throughout this paper and refer to our inverse problem as the ambient potential identification problem. We may drop the adjective ``ambient'' when there is no ambiguity. 

At first glance, the problem may appear less challenging since the ambient potential is independent of the population and remains invariant over time. 
However, the Nash equilibrium exhibits a highly nonlinear dependence on the ambient potential, and the well-posedness of the ambient potential identification problem remains an open question.
Although several numerical methods have been proposed for ambient potential identification, most of them reformulate the problem as a regularized optimization task and solve it using standard optimization algorithms. The underlying mechanism of the inverse problem remains poorly understood.

\paragraph{Contributions} 
We propose iterative algorithms for the ambient potential identification problem. These algorithms differ from traditional optimization-based methods and provide a new perspective on inverse mean-field games.
Our numerical experiments show that the Nash equilibrium measurement $\phi_0$ retains sufficient information to enable accurate reconstruction of the ambient potential, even when the viscosity $\nu$ is large. Our algorithms effectively leverage this information, achieving high reconstruction accuracy with relative errors down to $10^{-9}$ in a few iterations.
We further emphasize the connection between ambient potential identification in MFGs and the classical potential identification problem for linear parabolic equations. The Equilibrium Correction Iteration (ECI) Algorithm is closely related to fixed-point iteration methods developed for such inverse problems. This connection provides valuable insight into the structure of our approach and suggests promising directions for establishing well-posedness and convergence in future work.

\paragraph{Related Work}

The study of inverse problems in mean-field games has gained increasing attention in recent years. 

The theoretical advancements in this area primarily follow two main directions and are reviewed in \cite{liu2025inverse}.
One line of research extends the cost function to complex domains and employs harmonic sequences to study well-posedness. This approach was initially proposed in \cite{liu2022inverse} for mean field games (MFGs) with measurement $\phi_0$, and was subsequently extended to various other configurations in related works. A recent enhancement of this methodology is presented in \cite{liu2025decoding}.
A common and central assumption in these studies is that the running cost satisfies $q(x)+f(x,\rho)=0$ whenever $\rho=0$ for any $x,t$. However, this assumption limits their applicability to the ambient potential identification problem. 
To resolve the limitation, \cite{ren2024note} proves the well-posedness of the problem when $q$ is a constant, independent of both the state $x$ and time $t$. While insightful, this assumption is overly restrictive for practical applications.
The second direction leverages Carleman estimates to study the stability and uniqueness of inverse mean-field games. Following the framework introduced in \cite{klibanov2023lipschitz}, a series of works explore these properties in greater depth. 
The thesis \cite{liu2023inverse} investigates recovering the Lagrangian from measurements on the space boundary.
In addition to the above directions, \cite{yu2024invmfg_blo} proves the well-posedness of the ambient potential identification problem with a different measurement and from the perspective of numerical stability.

From the numerical perspective, most existing works formulate the MFG inverse problem as an optimization problem and solve it using optimization algorithms. 
For instance, \cite{ding2022mean,chow2022numerical} and \cite{yu2024invmfg_blo} address inverse MFGs using measurements of the equilibrium density and strategy $(\rho,v)$. The first two studies reformulate the inverse problem as a PDE-constrained optimization problem and solve it using primal-dual methods, while the latter focuses on potential MFGs, framing the inverse problem as a bilevel optimization problem and solving it with a gradient-based bilevel solver. 
Other approaches include the Gaussian Process framework proposed in \cite{guo2024decoding} and accelerated in \cite{yang2025gaussian}.
Additionally, \cite{klibanov2023convexification} aims to identify the interaction coefficient by expressing the unknown coefficient in terms of the Nash equilibrium and the observation, embedding it into the forward problem, and then using optimization algorithms to minimize the system's residual.
The only exception we are aware of is \cite{ren2024note}, which adopts the policy iteration forward solver developed in \cite{cacase2021policyitercvg}, and proposes a policy iteration inverse solver that alternates between solving an inverse linear parabolic equation and a forward FP equation.
While the implementation still relies on an optimization solver to solve the inverse linear parabolic equation, \cite{ren2024policy} provides a theoretical guarantee by proving linear convergence of the policy iteration inverse solver for the ambient potential identification problem.

In addition to the works discussed above, the ambient potential identification problem is closely related to the potential identification problem for linear parabolic equations. 
This is because when the viscosity is positive and the Hamiltonian is quadratic, the MFG system can be transformed into a system of semilinear parabolic equations via the Hopf-Cole transformation \cite{gueant2011mean,gueant2012hopfcole}. 
Notably, up to a Hopf-Cole transformation, our ECI algorithm has a similar formulation to fixed-point iteration methods \cite{Isakov2017bookinvprobforpde,isakov1991inverse,zhang2022identification} developed for potential identification in linear parabolic equations.

\paragraph{Organization} The rest of the paper is organized as follows. 
In Section \ref{sec: prob form}, we review the fundamental concepts of mean-field games and formulate the ambient potential identification problem. 
Section \ref{sec: algo} proposes the iterative algorithms for solving the inverse problem after reviewing the forward solver fictitious play. 
Section \ref{sec: env pot id unique} discusses the relationship between ambient potential identification in MFGs and inverse linear parabolic equations, and highlights how our ECI algorithm relates to fixed-point methods used for such inverse problems.
In Section \ref{sec: disc}, we describe the discretization and implementation details of the proposed methods. 
Section \ref{sec: num result} presents numerical results to demonstrate the effectiveness of our algorithms.
Finally, Section \ref{sec: conclusion} concludes this work and discusses future work.

\section{Ambient Potential Identification in Mean-Field Games}
\label{sec: prob form}

This section first reviews different formulations of the best response, defines the MFNE as the fixed point of the best response and formulates the forward problem. 
Then we formulate the ambient potential identification problem in mean-field games.

\subsection{The Best Response and Mean-Field Games}
\label{subsec: best response and mfg}
Consider an MFG system defined spatially on the $ d$-dimensional flat torus $\bbT^d$ and temporally on the interval $[0,T],(T>0)$. Let $\calP(\bbT^d)$ be the set of probability distributions on $\bbT^d$.
For any $t \in [0,T]$, $\rho(\cdot,t) \in \calP(\bbT^d)$ represents the state distribution at time $t$.
We slightly abuse our notation and use $\rho$ to represent a probability distribution or its density function, depending on the context.

For a given \textbf{population} state density flow $\tilde{\rho}\in C([0,T];\calP(\bbT^d))$, a \textbf{representative} player aims to solve an optimal control problem where the running cost and terminal cost depend on the population
\begin{equation}
\begin{aligned}
    \inf_{v} \quad & \mathbb{E}\left[ \int_0^T 
    \tilde{L}(X_s,v_s,\tilde{\rho}_s)\dd s + f_T(X_T,\tilde{\rho}_T)\dd x \right] \\
    &\text{subject to } X_0\sim \rho_0,\quad \dd X_s = v(X_s,s)\dd s + \sqrt{2\nu}\dd W_s.
\end{aligned}
\label{eq: individual oc}
\end{equation}
Here $\tilde{L}:\bbT^d \times \bbR^d \times \calP(\bbT^d) \to\bbR$ is the running cost assuming convex in $v$,
and $f_T:\bbT^d\times\calP(\bbT^d)\to\bbR$ is the terminal cost. 
$\rho_0$ is a given initial distribution, $\nu>0$ and $W_t$ is the standard Wiener process.
The optimal control problem \eqref{eq: individual oc} is equivalent to 
\begin{equation}
\begin{aligned}
    \inf_{\rho,v} \quad & \int_0^T \int_{\bbT^d}
    \rho_s(x)\tilde{L}(x,v_s,\tilde{\rho}_s)\dd x\dd s + \int_{\bbT^d} \rho_T(x) f_T(x,\tilde{\rho}_T) \dd x \\
    &\text{subject to } \rho(\cdot,0) = \rho_0,\quad \partial_t\rho -\nu\Delta\rho + \nabla\cdot(\rho v) = 0.
\end{aligned}
\label{eq: individual oc pdecstr}
\end{equation}
Here $\rho_t$ is the individual state distribution and $\tilde{\rho}_t$ is the population state distribution. The minimizer $(\rho,v)$ are the \textbf{best response} to the population distribution $\tilde{\rho}$.

By Pontryagin's maximum principle, the best response strategy is 
\begin{equation*}
    v(x,t) = - D_p\tilde{H}(x,\nabla\phi(x,t),\tilde{\rho}_t).
\end{equation*}
where $\phi$ is the value function and solves the Hamilton-Jacobi-Bellman (HJB) equation involving $\tilde{\rho}$,
\begin{equation}\left\{
    \begin{aligned}
        -\partial_t\phi(x,t)-\nu\Delta\phi(x,t)+\tilde{H}(x,\nabla\phi(x,t),\tilde{\rho}(\cdot,t))=0,
        &\quad (x,t)\in \bbT^d\times[0,T],  \\
        \phi(x,T)=f_T(x,\tilde{\rho}(\cdot,T)),
        &\quad x\in\bbT^d.    
    \end{aligned}\right.
    \label{eq: HJB}
\end{equation}
The best response distribution $\rho$ is induced by $v$ and therefore solves the Fokker-Planck (FP) equation
\begin{equation}\left\{
    \begin{aligned}
        \partial_t\rho(x,t)-\nu\Delta\rho(x,t)-\nabla\cdot(\rho D_p\tilde{H}(x,\nabla\phi,\tilde{\rho})(x,t))=0,
        \quad &(x,t)\in \bbT^d\times[0,T],  \\
        \rho(x,0)=\rho_0(x), 
        \quad &x\in\bbT^d.    
    \end{aligned}\right.
    \label{eq: FP}
    \end{equation}
If all individual players follow the best response strategy, the population distribution becomes the best response distribution $\rho$.
By definition, Nash Equilibrium is the fixed-point of the best response, i.e., the population distribution $\rho$ is the best response to itself.
Therefore, it can be formulated as the coupled PDE system \eqref{eq: mfg general}.
For the existence and uniqueness of the solution to the MFG system, we refer to \cite{Achdou2021MFGbook} and the references therein.

\subsection{The Ambient Potential Identification Problem}
\label{subsec: obs id prob}

In this paper, we restrict our focus to a separable cost 
\begin{equation*}
    \tilde{L}(x,v,\rho)=L(v)+q(x)+f(x,\rho).
\end{equation*}
$L$ is the dynamic cost, $f$ is referred to as the interaction cost because it represents the individual cost depending on the population, while $q$ is called the ambient potential, determined only by the environment.
Specifically, $q(x)$ quantifies the cost incurred to an individual at state $x$. 
The separable cost leads to the following MFG system 
\begin{equation}
\left\{\begin{aligned}
    -\partial_t\phi(x,t)-\nu\Delta\phi(x,t)+H(\nabla\phi)=q(x)+f(x,\rho(\cdot,t)),
    \quad &(x,t)\in \bbT^d\times[0,T],  \\
    \partial_t\rho(x,t)-\nu\Delta\rho(x,t)-\nabla\cdot(\rho D_pH(\nabla\phi))(x,t)=0,
    \quad &(x,t)\in \bbT^d\times[0,T],  \\
    \phi(x,T)=f_T(x,\rho(\cdot,T)),\quad \rho(x,0)=\rho_0(x), 
    \quad &x\in\bbT^d.    
\end{aligned}\right.
\label{eq: mfg obs}
\end{equation}

Assuming that the initial density $\rho_0$, interaction cost $f$ and terminal cost $f_T$ are known, we summarize the forward problem and inverse problem as follows.
\begin{problem}[Mean-Field Games with Ambient Potential]
\label{prob: forward}
    For given $q$, $(\rho,\phi;q)$ solves the forward problem if $(\rho,\phi)$ solves \eqref{eq: mfg obs} with ambient potential $q$.
\end{problem}
A constant shift in the ambient potential $q$ affects the forward solution $(\rho,\phi)$ as follows:
\begin{proposition}
    Let $c\in\mathbb{R}$ be a constant and $\phi_c(x,t) = \phi(x,t) + c(T-t)$, 
    then $(\rho,\phi;q)$ solves the forward problem if and only if $(\rho,\phi_c;q+c)$ solves the forward problem.
\label{prop: forward cst}
\end{proposition}

Given the potential $q$ and the interaction cost $f$, the system \eqref{eq: mfg obs} can be solved to determine the Nash equilibrium $(\rho, \phi)$. However, in many applications, the potential $q$ cannot be measured directly, even though some information about the Nash equilibrium may be observable. This motivates the consideration of the inverse problem: identifying the potential $q$ from measurements related to the Nash equilibrium.

\begin{problem}[Ambient Potential Identification in Mean-Field Games]
\label{prob: inverse}
    For given ${\phi}_0:={\phi}(\cdot,0)$,
    $(\hat{\rho},\hat{\phi},\hat{q};{\phi}_0)$ solves the inverse problem if $\hat{\phi}_0=\phi_0$ and $(\hat{\rho},\hat{\phi};\hat{q})$ solves the forward problem.
\end{problem}

\begin{remark}
    \cite{yu2024invmfg_blo} also investigates the ambient potential identification in MFGs, but their approach is based on measuring the density and velocity field $(\rho_t,v_t)$ at any $t\in[0,T]$.
    While $\rho_t,v_t$ are observable, requiring the information at all times $t$ may be too strong of an assumption.
    In contrast, our approach only requires the value function $\phi$ at $t=0$, which has the same dimension as the potential $q$.
\end{remark}

Similar to Proposition \ref{prop: forward cst}, it is straightforward to see how a constant shift in the potential $q$ affects the inverse problem.
\begin{proposition}
    $(\hat{\rho},\hat{\phi},\hat{q};\phi_0)$ solves the inverse problem if and only if $(\hat{\rho},\hat{\phi}_c,\hat{q}+c;\phi_0+cT)$ solves the inverse problem,
    where $c\in\mathbb{R}$ is a constant and $\phi_c(x,t) = \phi(x,t) + c(T-t)$.
\label{prop: inverse cst}
\end{proposition}

\section{Iterative Algorithms for Environment Potential Identification}
\label{sec: algo}

This section introduces iterative algorithms for the Ambient Potential Identification Problem (Problem \ref{prob: inverse}). 

We begin by reviewing the fictitious play algorithm \cite{cardaliaguet2017ficplay,yu2024ficplay}, an efficient forward solver for mean-field games based on best response dynamics. In fictitious play, each iteration computes the best response to the current population distribution and updates the population accordingly to approach the Nash equilibrium.
Building on this idea, we propose the Equilibrium Correction Iteration (ECI, Algorithm \ref{alg: ECI}). In ECI, each iteration solves the forward MFG system for a given estimate of the ambient potential and then updates this estimate by comparing the resulting Nash equilibrium with the observed measurement.
To further improve computational efficiency, we introduce two variants: the Best Response Iteration (BRI), which incorporates an inexact forward solver, and the Hierarchical Equilibrium Correction Iteration (HECI), which leverages a hierarchy of spatial grids.

\subsection{Fictitious Play through the Lens of Best Response}
\label{subsec: fictitious play}
Fictitious play \cite{cardaliaguet2017ficplay,yu2024ficplay} is an efficient iterative algorithm for solving the forward MFG by decoupling the FP and HJB equations. 
Each iteration of fictitious play consists of a best response step and an averaging step:
\begin{equation}
\begin{aligned}
    \text{Best response step:} \quad
    &(\rho^{(n)},v^{(n)}):=\argmin_{(\rho,v)\in\mathcal{C}} \mathcal{L}(\rho,v;\tilde{\rho}^{(n)})\\
    \text{Averaging step:} \quad
    &\tilde{\rho}^{(n+1)}=\left(1-\delta_n\right)\tilde{\rho}^{(n)}+\delta_n{\rho}^{(n)},
\end{aligned}
\label{eq: ficplay}
\end{equation}
where $\mathcal{C}$ is the set of admissible controls
\begin{equation*}
    \mathcal{C}= \{ (\rho,v):\rho(\cdot,0)=\rho_0,\partial_t\rho -\nu\Delta\rho + \nabla\cdot(\rho v) = 0 \},
\end{equation*} 
$\mathcal{L}$ is the objective of individual optimal control depending on population distribution 
\begin{equation*}
    \mathcal{L}(\rho,v;\tilde{\rho}) = \int_0^T\int_{\bbT^d}
    \rho_s(x) (L(v_s(x)) + q(x) + f(x,\tilde{\rho}_s) )\dd x\dd s 
    +\int_{\bbT^d} \rho_T(x)f_T(x,\tilde{\rho}_T)\dd x,
\end{equation*}
and $\delta_n$ is a weight for the averaging step.
As discussed in Section \ref{subsec: best response and mfg}, solving for the best response only requires solving the HJB and FP equations separately.
Therefore, the fictitious play algorithm decouples the system.
We summarize it in Algorithm \ref{alg: ficplay}.

As highlighted in \cite{yu2024ficplay}, the averaging step \eqref{eq: ficplay avg} and the choice of weight $\delta_n$ are crucial for the convergence of fictitious play. The weight $\delta_n$ can be preselected as a fixed value $\delta \in (0,1)$, set to $\frac{2}{n+2}$, or determined adaptively using backtracking line search.
When the interaction cost $f$ and the terminal cost $f_T$ are monotone, a sufficiently small weight $\delta_n$ guarantees that the updated state $\tilde{\rho}^{(n+1)}$ moves closer to the Nash equilibrium. Within a certain range determined by the problem structure and the current iterate $\rho^{(n)}$, smaller values of $\delta_n$ lead to more conservative updates and slower convergence, while larger values result in greater improvements and faster convergence. However, if $\delta_n$ is chosen too large, improvement toward the Nash equilibrium is no longer guaranteed, and the algorithm may diverge.
Under these monotonicity assumptions, Algorithm~\ref{alg: ficplay} with $\delta_n = \frac{2}{n+2}$ achieves sublinear convergence, as $\delta_n$ decays asymptotically. On the other hand, using a constant weight $\delta_n$ within a suitable range can lead to linear convergence, and adaptive choices of $\delta_n$ are designed to maximize progress at each iteration.
For further details, we refer to \cite{yu2024ficplay}.

\begin{algorithm}[ht]
\caption{Fictitious Play}
\begin{algorithmic}
    \STATE{Parameters} $\rho_0,f,f_T,\nu,0 < \delta_k \leq 1$
    \STATE{Initialization} $\tilde{\rho}^{(0)}$
    \FOR{$n=0,1,2,\cdots,N-1$}
    
        \STATE{Solve the HJB \eqref{eq: ficplay hjb} for $\phi^{(n)}$}
        \begin{equation}\left\{
        \begin{aligned}
            &-\partial_t\phi^{(n)}-\nu\Delta\phi^{(n)}+H(\nabla\phi^{(n)}(x,t))=q(x)+f(\tilde{\rho}^{(n)}(t)), \\
            & \phi^{(n)}(x,T)=f_T(x,\tilde{\rho}^{(n)}(T)).   
        \end{aligned}\right.
        \label{eq: ficplay hjb}
        \end{equation}

        \STATE{Solve the FP \eqref{eq: ficplay fp} for $\rho^{(n)}$}
        \begin{equation}\left\{
        \begin{aligned}
            &\partial_t\rho^{(n)}-\nu\Delta\rho^{(n)}-\nabla\cdot(\rho^{(n)} D_pH(x,\nabla\phi^{(n)}))=0, \\
            &\rho^{(n)}(x,0)=\rho_0(x).    
        \end{aligned}\right.
        \label{eq: ficplay fp}
        \end{equation}

        \STATE{Execute density average to obtain $\tilde{\rho}^{(n+1)}$}
        \begin{equation}
            \tilde{\rho}^{(n+1)}=\left(1-\delta_n\right)\tilde{\rho}^{(n)}+\delta_n\rho^{(n)}.
        \label{eq: ficplay avg}
        \end{equation}
    \ENDFOR
    \STATE{Output} $\tilde{\rho}^{(N)},\phi^{(N-1)}.$
\end{algorithmic}
\label{alg: ficplay}
\end{algorithm}

\subsection{The Equilibrium Correction Iteration}

Notice that from the HJB equation at $t=0$, 
\begin{equation}
    q = -\partial_t\phi|_{t=0} - \nu\Delta\phi_0 + H(\nabla\phi_0) - f(\rho_0),
\end{equation}
where $\rho_0$ is known and $\phi_0$ is given as the measurement.
If $\partial_t\phi|_{t=0}$ is also measurable, then from the HJB equation, we can directly obtain $q$.
While $\partial_t\phi|_{t=0}$ is unknown for the true potential $q$, we can solve the forward problem \eqref{eq: mfg obs} with some estimation of the environment potential, say $\hat{q}$, to obtain the solution $\hat{\phi}$ and then update the estimation by
\begin{equation}
    \hat{q}^+=-\partial_t\hat{\phi}|_{t=0}-\nu\Delta\phi_0+H(\nabla\phi_0)-f(\rho_0).
\label{eq: update subst dtphi term}
\end{equation}
Here $-\partial_t\hat{\phi}|_{t=0}$ is from the estimation and $-\nu\Delta\phi_0+H(\nabla\phi_0)$ are from the measurement.
Since $\hat{\phi}$ solves the forward problem \eqref{eq: mfg obs} with $\hat{q}$, this update is equivalent to
\begin{equation}
    \hat{q}^+=\hat{q}+(-\nu\Delta\phi_0+H(\nabla\phi_0))-(-\nu\Delta\hat{\phi}_0+H(\nabla\hat{\phi}_0)).
\label{eq: update subst dtphitilde term}
\end{equation}
Proposition \ref{prop: forward cst} reveals that $\nabla\phi$ and $\Delta\phi$ do not reflect a constant change in the ambient potential.
To incorporate the constant, we add an alignment term and update $\hat{q}$ with
\begin{equation}
    \hat{q}^+=\hat{q}+\frac{1}{T}\int_{\bbT^d}(\phi_0-\hat{\phi}_0)\dd x+(-\nu\Delta\phi_0+H(\nabla\phi_0))-(-\nu\Delta\hat{\phi}_0+H(\nabla\hat{\phi}_0)).
\label{eq: update incorp cst}
\end{equation}
We summarize the procedure in Algorithm \ref{alg: ECI} and name it the Equilibrium Correction Iteration (ECI).

\begin{algorithm}[ht]
\caption{The Equilibrium Correction Iteration (ECI)}
\begin{algorithmic}
    \STATE{Parameters} $\rho_0,f,f_T,\nu,\phi_0$
    \STATE{Initialization} $q^{(0)}$
    \FOR{$k=0,1,2,\cdots,K-1$}
    
        \STATE{update $\rho^{(k)},\phi^{(k)}$ by solving \eqref{eq: mfg obs} with potential $q^{(k)}$
        \begin{equation}\left\{
        \begin{aligned}
            &-\partial_t\phi^{(k)}(x,t) -\nu\Delta\phi^{(k)}(x,t) + H(\nabla\phi^{(k)}) = q^{(k)}(x) + f(\rho_t^{(k)}),\\
            &\partial_t\rho^{(k)}(x,t) -\nu\Delta\rho^{(k)}(x,t) -\nabla\cdot(\rho^{(k)}D_pH(\nabla\phi^{(k)})) = 0,\\
            &\phi^{(k)}(x,T)=f_T(x,\rho_T^{(k)}), \quad \rho^{(k)}(\cdot,0)=\rho_0.
        \end{aligned}\right.
        \label{eq: eci mfg}
        \end{equation}
        }

        \STATE{update $q^{(k+1)}$ with
        \begin{equation}
        \begin{aligned}
            q^{(k+1)} = q^{(k)} 
            &+ \left( \frac{1}{T}\int_{\bbT^d}\phi_0(x)\dd x-\nu\Delta\phi_0 + H(\nabla\phi_0)\right) \\
            &- \left( \frac{1}{T}\int_{\bbT^d}\phi_0^{(k)}(x)\dd x -\nu\Delta\phi_0^{(k)} + H(\nabla\phi_0^{(k)})\right)
        \end{aligned}
        \end{equation}
        
        }
        
    \ENDFOR
    \STATE{Output} $q^{(K)},\rho^{(K-1)},\phi^{(K-1)}.$

\end{algorithmic}
\label{alg: ECI}
\end{algorithm}

The following proposition is straightforward.

\begin{proposition}
\label{prop: sol -> fixed point}
If $(\hat{\rho},\hat{\phi},\hat{q};{\phi}_0)$ solves the inverse problem and $(\hat{\rho},\hat{\phi})$ uniquely solves the forward problem with $\hat{q}$, then $\hat{q}$ is a fixed point of Algorithm \ref{alg: ECI}.
\end{proposition}
At this moment, we do not have proof of the reverse proposition or the convergence of the algorithm under reasonable assumptions.
But we provide some reasoning here and support them with numerical experiments in Section \ref{subsec: num cvg}.
We define the correction of the estimation $\hat{q}$ to ground truth $q$ as
\begin{equation}
    c(q,\hat{q}) := (q - \hat{q}) + \frac{1}{T}\int_{\bbT^d}(\phi_0-\hat{\phi}_0)\dd x + \partial_t(\phi - \hat{\phi})\big|_{t=0}.
\label{eq: update correction term}
\end{equation}
Then the update \eqref{eq: update incorp cst} is equivalent to:
\begin{equation}
    \hat{q}^+ = \hat{q} + c(q,\hat{q}).
\label{eq: update in terms of correction}
\end{equation}
Notice that if the terms involving $\phi$ and $\hat{\phi}$ are absent in the correction, then the update is exact.
Therefore, we define the error of the correction as
\begin{equation}
    e(\phi,\hat{\phi}) := \frac{1}{T}\int_{\bbT^d}(\phi_0-\hat{\phi}_0)\dd x + \partial_t(\phi - \hat{\phi})\big|_{t=0}.
\label{eq: update error term}
\end{equation}
The correction term $c(q,\hat{q})$ and the error in correction $e(\phi,\hat{\phi})$ provide useful diagnostics for the behavior of a single ECI update. 
Note that the error in correction is exactly the error of the updated estimator $\hat{q}^+$:
\begin{equation}
    \hat{q}^+ - q = e(\phi,\hat{\phi}).
\label{eq: error of correction is error of new est}
\end{equation}
Thus, if $\|e(\phi,\hat{\phi})\|$ is sufficiently small, the new estimator $\hat{q}^+$ is a sufficiently accurate estimator of the true potential $q$.
If $\|e(\phi,\hat{\phi})\|$ is not negligible but satisfies
\begin{equation}
    \|e(\phi,\hat{\phi})\| < \|q - \hat{q}\|,
\end{equation}
then the ECI update reduces the estimation error.
If at some state $x$,
\begin{equation}
    (q - \hat{q}) (q - \hat{q} + e(\phi,\hat{\phi})) \geq 0,
\label{eq: eci good cond}
\end{equation}
the correction term is aligned with the error and pushes the estimator in the correct direction.

These indicators are examined in the numerical experiments of Section~\ref{subsec: num cvg}, which show that the ECI algorithm effectively recovers the potential $q$, particularly when the interaction cost $f$ is monotone, the viscosity $\nu$ is large, and the terminal time $T$ is small.

We admit that the ECI algorithm relies on $\nabla \phi_0$ and $\Delta \phi_0$, and is therefore sensitive to noise in the observations of $\phi_0$. However, numerical experiments in Section \ref{sec: num result}, along with the connection to the inverse linear parabolic equation discussed in Section \ref{sec: env pot id unique}, suggest that the information of $q$ is accurately encoded in $\phi_0$ through $H(\nabla \phi_0)$ and $\Delta \phi_0$, and that ECI effectively leverages this information.
Whether this information is significantly degraded by noise, and if not, how to robustly extract it, remain important and promising directions for future research.

\subsection{The Best Response Iteration}
\label{subsec: bri}

Each iteration of the ECI method involves solving a forward MFG system \eqref{eq: eci mfg}, which is computationally challenging due to the interdependence and forward-backward structure of the HJB and FP equations. 
To address this, we integrate the forward fictitious play solver (Algorithm \ref{alg: ficplay}) into ECI (Algorithm \ref{alg: ECI}) and propose the Best Response Iteration (BRI) in Algorithm \ref{alg: BRI}. 
To be precise, for the given ambient potential estimate, instead of solving for the exact Nash Equilibrium, we only run the forward solver for a fixed number of iterations to obtain an approximate Nash Equilibrium.

BRI is a double-loop algorithm: the inner loop output $\phi^{(k,N)}$ approximates the forward problem solution $\phi^{(k)}$ for the estimated potential $q^{(k)}$, while the outer loop updates the ambient potential based on the approximate estimated measurement. ECI can be viewed as a special case of BRI where $N=\infty$. Numerical experiments in Section \ref{subsec: num comp} indicate that BRI remains effective even with $N=1$. This behavior is consistent with other double-loop algorithms, such as the bilevel solver in \cite{yu2024invmfg_blo}. The key mechanism is that the accuracy of both the inner and outer loops improves as iterations progress. 
As a result, BRI significantly reduces the computational cost of ECI while maintaining its effectiveness.

We emphasize that the inner loop is designed to approximate the solution of the forward problem, and the performance of BRI depends on how accurately the inner loop approximates this solution. For the BRI method with a fixed number of inner iterations $N$ to serve as an effective acceleration of ECI, it is crucial to initialize the inner loop with $\tilde{\rho}^{(k,0)} = \rho^{(k-1)}$ and to choose $\delta_{k,n}$ properly.
Since $\rho^{(k-1)}$ approximates the Nash equilibrium associated with $q^{(k-1)}$, and $q^{(k)}$ is updated from $q^{(k-1)}$, initializing the inner loop with $\rho^{(k-1)}$ and running a fixed number of fictitious play steps provides a better approximation of the Nash equilibrium corresponding to $q^{(k)}$. In contrast, initializing $\tilde{\rho}^{(k,0)}$ with a fixed distribution $\rho$, regardless of the iteration, may prevent the inner loop from effectively reducing the forward residual. Consequently, BRI may converge to an incorrect potential.
Numerical experiments in Section \ref{subsec: num bri} demonstrate that using $\tilde{\rho}^{(k,0)} = \rho^{(k-1)}$ is necessary for BRI to converge to the correct potential.
Similarly, the choice of $\delta_{k,n}$ also affects BRI, as it influences both the convergence and convergence rate of the forward solver, fictitious play. Therefore, it affects how accurately $\tilde{\rho}^{(k,N)}$ approximates the Nash equilibrium, which in turn impacts the convergence of BRI.
In our numerical experiments in Section~\ref{subsec: num comp}, we show that when $N = 1$, using $\delta_{k,n} = 0.2$ yields smaller improvements in the forward residual and therefore slower convergence of BRI compared to $\delta_{k,n} = 0.5$. 
However, this does not imply that $\delta_{k,n}$ should be chosen arbitrarily large. As discussed in the convergence analysis of \cite{yu2024ficplay}, an excessively large $\delta_{k,n}$ may fail to reduce the forward residual and thus hinder potential identification. Numerical results in Section~\ref{subsec: num bri} further show that too large a value of $\delta_{k,n}$ can lead to oscillations in BRI.

\begin{algorithm}[htb]
\caption{The Best Response Iteration (BRI)}
\begin{algorithmic}
    \STATE{Parameters} $\rho_0,f,f_T,\nu,\phi_0$
    \STATE{Initialization} $q^{(0)},\rho^{(-1)}$
    \FOR{$k=0,1,2,\cdots,K-1$}
        \STATE{Assign $\tilde{\rho}^{(k,0)}=\rho^{(k-1)}$}
        \FOR{$n=0,1,2,\cdots,N-1$}
        
            \STATE{Solve the HJB for $\phi^{(k,n)}$ 
            \begin{equation}\left\{
            \begin{aligned}
                &-\partial_t\phi^{(k,n)}(x,t) 
                 -\nu\Delta\phi^{(k,n)}(x,t) 
                 + H(\nabla\phi^{(k,n)}(x,t)) 
                 = q^{(k)}(x) + f(\tilde{\rho}_t^{(k,n)}),\\
                &\phi^{(k,n)}(x,T)=f_T(x,\tilde{\rho}_T^{(k,n)})
            \end{aligned}\right.
            \label{eq: BRI HJB}
            \end{equation}
            }
    
            \STATE{Solve the FP for ${\rho}^{(k,n)}$
            \begin{equation}\left\{
            \begin{aligned}
                &\partial_t{\rho}^{(k,n)}(x,t) 
                 -\nu\Delta{\rho}^{(k,n)}(x,t) 
                 -\nabla\cdot({\rho}^{(k,n)}D_pH(\nabla\phi^{(k,n)})) = 0,\\
                &{\rho}^{(k,n)}(\cdot,0)=\rho_0.\\
            \end{aligned}\right.
            \label{eq: BRI FP}
            \end{equation}        
            }

            \STATE{Update $\tilde{\rho}^{(k,n+1)}$ by
            \begin{equation}
                \tilde{\rho}^{(k,n+1)} = (1-\delta_{k,n})\tilde{\rho}^{(k,n)} + \delta_{k,n}\rho^{(k,n)}
            \label{eq: BRI avg}
            \end{equation}        
            
            }
            
        \ENDFOR
        \STATE{Assign $\rho^{(k)}=\tilde{\rho}^{(k,N)},\phi^{(k)}=\phi^{(k,N-1)}$}

        \STATE{update $q^{(k+1)}$ with
        \begin{equation}
        \begin{aligned}
            q^{(k+1)} = q^{(k)} 
            &+ \left( \frac{1}{T}\int_{\bbT^d}\phi_0(x)\dd x-\nu\Delta\phi_0 + H(\nabla\phi_0)\right) \\
            &- \left( \frac{1}{T}\int_{\bbT^d}\phi_0^{(k)}(x)\dd x -\nu\Delta\phi_0^{(k)} + H(\nabla\phi_0^{(k)})\right)
        \end{aligned}
        \label{eq: BRI update}
        \end{equation}
        }
        
    \ENDFOR
    \STATE{Output} $q^{(K)},\rho^{(K-1)},\phi^{(K-1)}.$

\end{algorithmic}
\label{alg: BRI}
\end{algorithm}

\subsection{Hierarchical ECI}
\label{subsec: hierarchical ECI}

Another approach to reduce the computational cost of ECI is to use a hierarchical grid as in \cite{yu2024ficplay}.
The idea is to solve the forward problem on a coarse grid and then interpolate the solution to a finer grid.
The solution on the coarse grid serves as a good initial guess for the fine grid. 
We summarize the hierarchical ECI (HECI) in Algorithm \ref{alg: HECI}.

As discussed in Section \ref{sec: disc}, we use a Lax-Friedrichs discretization which introduces a numerical viscosity that is proportional to the mesh step. 
Therefore, starting from the coarse grid can also stabilize the algorithm when the physical viscosity is small.
Numerical experiments in Section \ref{subsec: num acc} show that HECI is especially effective in capturing the low-frequency components of the ambient potential.

\begin{algorithm}[htb]
    \caption{The Hierarchical Equilibrium Correction Iteration (HECI)}
    \begin{algorithmic}
        \STATE{Parameters} $\rho_0,f,f_T,\nu,\phi_0$
        \STATE{Initialization} $q^{(0)},\rho^{(0)}$ on coarse grid $\calG_0$
        \STATE{Compute the observation-induced term on the desired fine grid $\calG_L$
        $$ \left( \frac{1}{T}\int_{\bbT^d}\phi_0(x)\dd x-\nu\Delta\phi_0 + H(\nabla\phi_0)\right)_{\calG_L} $$
        }
        \FOR{$l=1,\cdots,L$}
            \STATE{Restrict the observation-induced term to the current grid $\calG_l$}
            \STATE{Compute $q^{(l,0)}$ by interpolating $q^{(l-1)}$ to the current grid $\calG_l$}
            \STATE{Run ECI (Algorithm \ref{alg: ECI}) on the current grid $\calG_l$ with initialization $q^{(l,0)}$ and obtain output $ q^{(l,K)}$ and $\rho^{(l,K-1)},\phi^{(l,K-1)}$}
            \STATE{Denote the output of ECI on the current grid as $q^{(l)}:=q^{(l,K)}$}
        \ENDFOR
        \STATE{Output} $q^{(L)},\rho^{(L,K-1)},\phi^{(L,K-1)}.$
    
    \end{algorithmic}
    \label{alg: HECI}
    \end{algorithm}

\section{Connection to Inverse Linear Parabolic Equations}
\label{sec: env pot id unique}

This section explores the connection between our ambient potential identification problem and classical inverse problems for linear parabolic equations. We highlight how the structure of our ECI algorithm parallels fixed-point methods developed for such inverse problems. This perspective offers valuable insight into the underlying mechanism of our approach and suggests directions for future theoretical analysis, including well-posedness and convergence.

When $\nu>0$ and $H(p)=\frac{1}{2}|p|^2$, the MFG system \eqref{eq: mfg obs} can be transformed into a semilinear parabolic system using the Hopf-Cole transformation.  \cite{gueant2011mean,gueant2012hopfcole}.
\begin{proposition}[Hopf-Cole Transformation for MFGs \cite{gueant2011mean,gueant2012hopfcole}]
\label{prop: Hopf-Cole}
If $(\rho,\phi)$ solves the MFG system \eqref{eq: mfg obs} with quadratic Hamiltonian $H(p)=\frac{1}{2}|p|^2$ and $\nu>0$,
then 
$$(w,u)=\left(\exp{\left( -\frac{\phi}{2\nu} \right)}, \rho\exp{\left( \frac{\phi}{2\nu} \right)}\right)$$ 
solves the following system:
\begin{equation}\left\{
    \begin{aligned}
        -\partial_t w - \nu\Delta w + \frac{1}{2\nu}(q+f(wu))w=0, \quad &(x,t)\in\bbT^d\times[0,T],\\
        \partial_t u - \nu\Delta u + \frac{1}{2\nu}(q+f(wu))u=0,\quad &(x,t)\in\bbT^d\times[0,T],\\
        w(x,T) = \exp{\left( -\frac{f_T(x,(wu)_T)}{2\nu} \right)},\quad u(x,0) = \frac{\rho_0(x)}{w(x,0)},\quad &x\in\bbT^d.
    \end{aligned}\right.
    \label{eq: mfg linear parabolic}
    \end{equation}
\end{proposition}
\begin{proof}
    The boundary conditions are easy to verify.
    The equation for $w$ is straightforward by noticing
    $$ \partial_t w = -\frac{w}{2\nu}\partial_t\phi,\quad
    \nabla w = -\frac{w}{2\nu}\nabla\phi,\quad
    \Delta w = \frac{w}{2\nu^2}\left( -\nu\Delta\phi + \frac{1}{2}|\nabla\phi|^2 \right). $$
    Since $wu=\rho$, from the FP equation, we have
    $$ \partial_t (wu) - \nu\Delta(wu) + 2\nu\nabla\cdot(u\nabla w) = 0. $$
    Therefore, 
    $$ \partial_t u - \nu\Delta u - \frac{u}{w}(-\partial_t w - \nu\Delta w) = 0 $$
    and gives the equation for $u$.
\end{proof}

The Hopf-Cole transformation recasts the ambient potential identification problem as follows: given $u_0$, $f$, $f_T$, and the measurement $w_0$, recover $q$ and $(w,u)$. 

If $f_T$ is independent of $\rho$ but $f$ is present, we obtain initial and terminal data for $w$, which solves a semilinear parabolic equation. In this case, the standard theory and strategies for inverse linear parabolic equations does not directly apply because of the time-dependent term $f(wu)$. And the uniqueness, existence, and stability for the inverse problem \ref{prob: inverse} remain open questions.

If both $f_T$ and $f$ are independent of $\rho$ (i.e., $f$ is absent), the system decouples. Although this setting is restrictive and omits the main challenges of MFGs, it provides a useful starting point for understanding the ECI algorithm.
In this case, the ambient potential identification problem reduces to recovering $q$ from measurements $w_0$ and $w_T$, where $w$ solves
\begin{equation}
    -\partial_t w - \nu\Delta w + \frac{1}{2\nu}q w = 0.
\label{eq: linpara}
\end{equation}
This problem has been extensively studied; see, for example, 
\cite{isakov1991inverse} for results on existence, \cite{isakov1993uniqueness} for uniqueness and Chapter 9 of the \cite{Isakov2017bookinvprobforpde} for further related topics.
A constructive existence proof for $(w, q)$ is given in \cite{isakov1991inverse,Isakov2017bookinvprobforpde}, based on the monotonicity of the following iterative scheme:
\begin{equation}
    q^{(k+1)} = \frac{2\nu}{w_0}\left( \partial_t w^{(k)}\big|_{t=0} + \nu\Delta w_0 \right),
\label{eq: inv linpara update 1}
\end{equation}
where $w^{(k)}$ solves \eqref{eq: linpara} with $w(x,T)=w_T(x)$ and potential $q^{(k)}$. 
The idea is to use $\partial_t w^{(k)}|_{t=0}$, computed from the current estimate $q^{(k)}$, as an approximation to the true (unknown) $\partial_t w|_{t=0}$, and update $q$ accordingly. 
Using the equation for $w^{(k)}$, this update can also be written as
\begin{equation}
    q^{(k+1)} = \frac{w_0^{(k)}}{w_0}q^{(k)} + \frac{2\nu^2}{w_0}\Delta \left(w_0 - w_0^{(k)}\right),
\end{equation}
where $w_0^{(k)}$ denotes $w^{(k)}$ at $t=0$.
Applying the change of variables $w = \exp{\left( -\frac{\phi}{2\nu} \right)}$, the update can be rewritten in terms of the HJB measurement $\phi_0$ as
\begin{equation}
    q^{(k+1)} = r_0^{(k)}q^{(k)} + \left(-\nu\Delta\phi_0 + \frac{1}{2}|\nabla\phi_0|^2 \right) - r_0^{(k)}\left(-\nu\Delta\phi_0^{(k)} + \frac{1}{2}|\nabla\phi_0^{(k)}|^2 \right),
\label{eq: inv linpara update 2}
\end{equation}
where $r_0^{(k)} = \exp{\left( -\frac{\phi_0^{(k)}-\phi_0}{2\nu} \right)}$.

Comparing \eqref{eq: inv linpara update 2} with our ECI update \eqref{eq: update subst dtphitilde term}, we see they are structurally similar, except for the multiplicative factor $r_0^{(k)}$ in the linear parabolic case. 
While the monotonicity argument in \cite{Isakov2017bookinvprobforpde} can be adapted to analyze the ECI update, the required conditions are quite restrictive and may not hold in general nonlinear MFG settings. Establishing convergence of ECI under broader and more practical assumptions remains an interesting direction for future research.

\section{Discretization and Implementation}
\label{sec: disc}

We present the details of implementation in this section.
The settings are very close to those in \cite{yu2024ficplay}, except that the boundary conditions here are periodic.

Let $n_t,n_x$ be positive integers and $\dt=\frac{T}{n_t},t_n=n\dt, n=0,\cdots,n_t$, $\dx=\frac{1}{n_x},x_i=i\dx, i=0,\cdots,n_x$.
Denote $\calG$ as the collection of all grid points $(x_i,t_n),i=0,\cdots,n_x,n=0,\cdots,n_t$.
For any function $u$ defined on $\Omega\times[0,T]$, we denote the approximation of $u$ on $\calG$ as $u_{\calG}$ where $(u_{\calG})_\myij$ approximates $u(x_i,t_n)$ and omit $\calG$ in the subscript when the context provides no ambiguity.
In this paper, we consider the periodic boundary condition. To tackle this in the discretization, we set $u_0=u_{n_x}$, $u_{-1}=u_{n_x-1}$ and $u_{n_x+1}=u_{1}$. 
This results in the following finite difference operators on the grid $\calG$.
Firstly, the elementary one-sided finite difference operators are defined as
\begin{equation*}
    (\dxp u)_\myij = 
        \frac{u_\ipj-u_\myij}{\dx}, 
    \quad
    (\dxm u)_\myij = 
        \frac{u_\myij-u_\imj}{\dx},  
    \quad 
    i=0,\cdots,n_x.
\end{equation*}
Then we define the central difference operator $(\dxc u)_\myij$ as the average of the one-sided finite difference, and $[D_xu]_\myij$ as the collection of the one-sided finite differences at $(x_i,t_n)$
\begin{equation*}
    (\dxc u)_\myij = \frac{1}{2}(\dxp u)_\myij + \frac{1}{2}(\dxm u)_\myij,\quad
    (D_xu)_\myij = \left((\dxp u)_\myij,(\dxm u)_\myij\right),
\end{equation*}
and the discrete Laplace operator as
\begin{equation*}
    (\Delta_x u)_\myij =
        \frac{(\dxp u)_\myij-(\dxm u)_\myij}{\dx}, i=0,\cdots,n_x
\end{equation*}
We call $u=(u_\myij)_{i=0,n=0}^{i=n_x,n=n_t}$ a function on $\calG$ and $(v)=(v^+,v^-)$ a velocity field on $\calG$ where $v^+,v^-$ are functions on $\calG$.
For any functions $u_1,u_2$ and vector fields $(v_1)=(v_1^+,v_1^-),(v_2)=(v_2^+,v_2^-)$ on $\calG$, we define the inner product on the grid $\calG$ as
\begin{equation*}
    \langle u_1,u_2 \rangle_{\calG}:=\Delta t\Delta x\sum_{i=0}^{n_x}\sum_{n=0}^{n_t} (u_1)_\myij(u_2)_\myij,\quad
    \left\langle (v_1),(v_2) \right\rangle_{\calG}:=
    \frac{1}{2}\left(
    \left\langle v_1^+, v_2^+ \right\rangle_{\calG} + 
    \left\langle v_1^-, v_2^- \right\rangle_{\calG}
    \right)
\end{equation*}
The way we introduce the ghost point values makes the discrete Laplacian operator self-adjoint under the above inner product.
To preserve the adjoint relation between gradient and negative divergence, for a given velocity field $v$ on $\calG$, we define its divergence to be $-D_x^*[v]$ where $D_x^*$ is the adjoint operator of $D_x$,
\begin{equation*}
    (D_x^*(v))_\myij := 
        -\frac{1}{2}\left(\dxm(v^+) + \dxp(v^-)\right)_\myij, i=0,\cdots,n_x,
\end{equation*}
and then
\begin{equation*}
    \langle u,D_x^*(v) \rangle_{\calG} = \left\langle (D_x(u)),(v) \right\rangle_{\calG}.
\end{equation*}

We adopt the Lax-Friedrichs Hamiltonian, defined as:
\begin{equation*}
    \hamlf(x,\pxp,\pxm)
    =H\left( x, \left( \frac{\pxp+\pxm}{2} \right) \right)
    -\nu_n\frac{\pxp-\pxm}{2},
\end{equation*}
where $\nu_n\geq0$ is the numerical viscosity coefficient.

The interaction cost $f(x,\rho)$ and the terminal cost $f_T(x,\rho)$ at $x_i$ are approximated by $f_{\calG}(x_i,\rho_{\calG})$ and $f_{T,\calG}(x_i,\rho_{\calG})$, respectively. For any function $\rho_{\calG}$ defined on the grid $\calG$, we define its corresponding piecewise constant function on $\Omega$ as $\bar{\rho}_{\calG}$, where:
\[
\bar{\rho}_{\calG}(x) = \rho_i, \quad x_i-\frac{\dx}{2}\leq x < x_i+\frac{\dx}{2}.
\]
To ensure consistency, we assume there exists a constant $C$ independent of $\dx$ and $\rho$ such that for any $x_i$ and $\rho_{\calG}$:
\[
|f_{\calG}(x_i,\rho_{\calG})-f(x_i,\bar{\rho}_{\calG})|<C\dx,
\quad |f_{T,\calG}(x_i,\rho_{\calG})-f_T(x_i,\bar{\rho}_{\calG})|<C\dx.
\] 

Below are examples of cost functions and their consistent discretizations:
\begin{enumerate}
    \item If $f(x,\rho)=f_0(\rho(x))$ is defined locally by a function $f_0:\bbR\to\bbR$, we can approximate it as $f_{\calG}(x_i,\rho_{\calG})=f_0(\rho_i)$.
    \item If $f(x,\rho)=\int_\Omega K(x,y)\rho(y)\dd y$ is a nonlocal cost defined by a kernel function $K$, it can be approximated as $f_{\calG}(x_i,\rho_{\calG})=\sum_{j=0}^{n_x} K(x_i,x_j)\rho_j\dx$.
\end{enumerate}    

Using the above notation, we solve the discrete HJB equation as follows. First, we set the terminal condition:
\begin{equation}
    \phi_{i,n_t} = f_T(x_i,\rho_{\cdot,n_t}),\qquad i=0,1,\cdots,n_x.
\end{equation}
Then, for $\phi_{\cdot,n}$, where $n=n_t-1,n_t-2,\cdots,0$, we solve backward in time using Newton's iteration to address the nonlinear system:
\begin{equation}
    -\frac{\phi_\ijp-\phi_\myij}{\dt} - \nu(\Delta_x \phi)_\myij + \hamlf(x_i,(D_x\phi)_\myij)-f(x_i,\rho_{\cdot,n+1})=0,
    \quad i=0,1,\cdots,n_x.
\end{equation}

For the FP equation, we ensure the adjoint relationship between the HJB and FP equations by employing the following discretization:
\begin{equation}
    \frac{\rho_\ijp-\rho_\myij}{\dt} 
        - \nu(\Delta_x \rho)_\ijp + (D_x^*(\rho_{\cdot,n+1} v_{\cdot,n}))_i =0,
    \quad i=0,1,\cdots,n_x.
\end{equation}
Here, the velocity field is given by:
\[
v^+_\myij = -\nabla_{p^+}\hamlf(x_i,(D_x\phi)_\myij), \quad v^-_\myij = -\nabla_{p^-}\hamlf(x_i,(D_x\phi)_\myij).
\]
This corresponds to the Lax-Friedrichs scheme for the FP equation:
\[
\frac{\rho_\ijp-\rho_\myij}{\dt} 
- \left(\nu+\frac{\nu_n}{2}\dx\right)(\Delta_x \rho)_\ijp 
+ (D_x^*(\rho_{\cdot,n+1}{v}_{\cdot,n}))_i =0,
\]
where ${v}^{\pm}_\myij = -\nabla_pH(x_i,(\dxc \phi)_\myij)$. Notably, this scheme preserves discrete mass conservation.

\section{Numerical Results}
\label{sec: num result}

This section presents numerical results to verify the effectiveness of the proposed algorithms and the conjecture. All of our numerical experiments are implemented in MATLAB on a PC with an Apple M3 Pro chip and 18 GB of memory.

We denote the density function of a univariate Gaussian distribution as $$\rho_G(x;\mu,\sigma):=\frac{1}{\sqrt{2\pi}\sigma}\exp{\left(-\frac{(x-\mu)^2}{2\sigma^2}\right)}.$$
Let the inner product and norm on the discrete space grid $\calG$ be
$$\langle u,v \rangle_{\calG}:= \Delta x \sum_{i}u_{i}v_{i},\quad
\left\| u \right\|_{\calG} :=  \sqrt{\langle u,u \rangle_{\calG} }.$$
When there is no ambiguity, we omit the subscript $\calG$ in the notation.

We generate the observation data by solving the forward problem \eqref{eq: mfg obs} with a given potential $q$ with fictitious play or hierarchical fictitious play \cite{yu2024ficplay}.
We use the measurement relative error $\frac{\|\phi^{(k)}_0-\phi_0\|}{\|\phi_0\|}$ to assess the convergence of the algorithm and report the ambient potential relative error $\frac{\|q^{(k)} - q\|}{\|q\|}$ to evaluate the accuracy of the solution.
When comparing with the policy iteration method in \cite{ren2024policy}, we label it as PIT and the observation data for PIT are generated using the forward solver based on policy iteration \cite{cacase2021policyitercvg}.

\subsection{Effectiveness of ECI}
\label{subsec: num effective}

We first demonstrate that the measurement $\phi_0$ retains sufficient information to accurately reconstruct the ambient potential, even in the presence of significant viscosity. The ECI algorithm effectively leverages this information, allowing a rapid and highly accurate recovery of the true potential.

We consider the MFG on the spatial domain $[0,1]^d$ with periodic boundary conditions and $T=1$. 
For the 1D problem, $\Delta x=5\times 10^{-4},\Delta t=2\times10^{-3}.$ For the 2D problem, $\Delta t = 2\times10^{-2}$ and $\Delta x = \frac{1}{201},\Delta y=\frac{1}{170}$ to align with the image size $202\times 171$.
The Hamiltonian is $H(p)=\frac{1}{2}|p|^2$. The initial distribution is an isotropic Gaussian centered at $0.5$ (in 1D) or $(0.5,0.5)$ (in 2D) with variance $0.01$ in each direction. The interaction cost is $f(x,\rho)=K*\rho^2$ with $K$ being a Gaussian kernel. For the 1D example, $K(x,y)=\rho_G(x-y;0,0.1)$, and for the 2D example, the kernel is anisotropic $K(x,y)=\rho_G(x_1-y_1;0,0.1)\rho_G(x_2-y_2;0,0.2)$. The terminal cost is $f_T(x) = 0$.
The viscosity is set to $\nu=0.2$ for the 1D example and $\nu=0.1$ for the 2D example.

\begin{figure}[htb]
    \centering
        \begin{subfigure}[t]{0.48\textwidth}
            \centering
            \caption*{True potential and PIT recovered potential}
            \includegraphics[width=\textwidth]{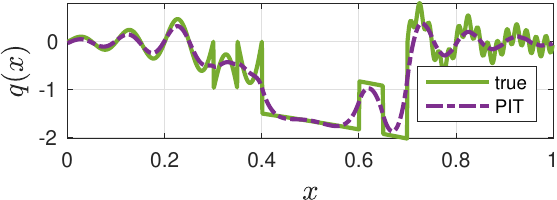}
        \end{subfigure}
        \begin{subfigure}[t]{0.48\textwidth}
            \centering
            \caption*{Measurement-induced term $M(\phi_0)$ in ECI}
            \includegraphics[width=\textwidth]{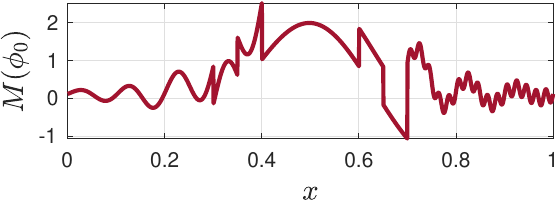}
        \end{subfigure}
        \begin{subfigure}[t]{0.48\textwidth}
            \centering
            \caption*{ECI result after 1 iter, rel. err. $2.41\times 10^{-1}$}
            \includegraphics[width=\textwidth]{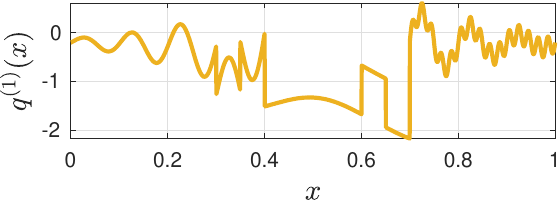}
        \end{subfigure}
        \begin{subfigure}[t]{0.48\textwidth}
            \centering
            \caption*{ECI result after 23 iters, rel. err. $2.88\times 10^{-9}$}
            \includegraphics[width=\textwidth]{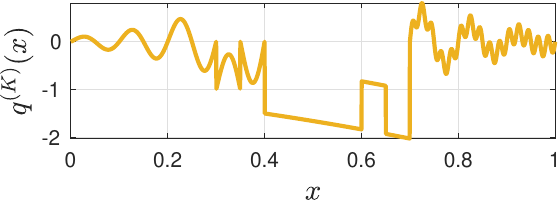}
        \end{subfigure}
        \begin{subfigure}[t]{0.48\textwidth}
            \centering
            \caption*{$q$ rel. err. vs num. of outer loop}
            \includegraphics[width=\textwidth]{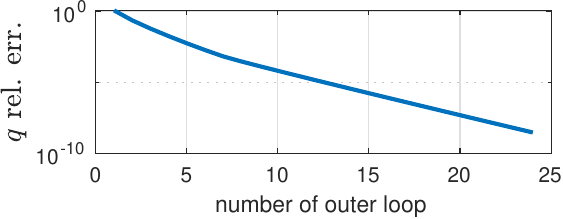}
        \end{subfigure}
        \begin{subfigure}[t]{0.48\textwidth}
            \centering
            \caption*{$\phi_0$ rel. err. vs num. of outer loop}
            \includegraphics[width=\textwidth]{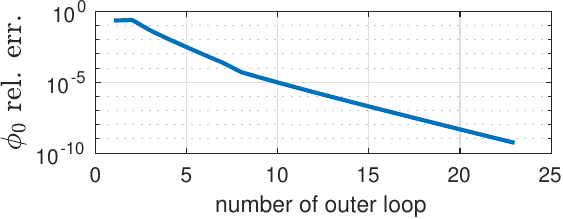}
        \end{subfigure}
    \caption{ECI efficacy with 1D example (\S\ref{subsec: num effective}).  
    The measurement-induced term $M(\phi_0):=\frac{1}{T}\int_{\mathbb{T}^d}\phi_0(x)\dd x-\nu\Delta\phi_0+H(\nabla\phi_0)$ retains substantial information about the true potential $q$. 
    The ECI algorithm utilizing $M(\phi_0)$ rapidly recovers the main structure of the true potential. 
    }
    \label{fig: 1dillu}
\end{figure}

\begin{figure}[htb]
    \centering
        \begin{subfigure}[t]{0.32\textwidth}
            \flushleft
            \caption*{True $q$}
            \includegraphics[width=\textwidth]{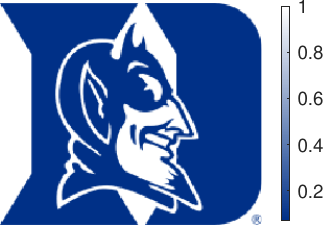}
        \end{subfigure}
        \begin{subfigure}[t]{0.32\textwidth}
            \flushleft
            \caption*{meas.-induced term $M(\phi_0)$ in ECI}
            \includegraphics[width=\textwidth]{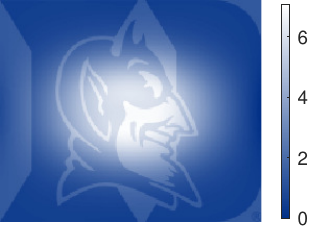}
        \end{subfigure}
        \begin{subfigure}[t]{0.32\textwidth}
            \flushleft
            \caption*{$q$ rel. err. vs num. outer it.}
            \includegraphics[width=\textwidth]{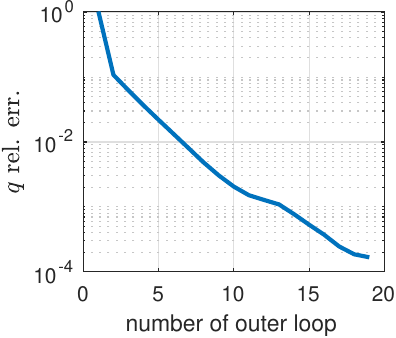}
        \end{subfigure}
    \vspace{1em}
        \begin{subfigure}[t]{0.32\textwidth}
            \flushleft
            \caption*{After 1 iter $q^{(1)}$, rel. err. $1.08\times 10^{-1}$}
            \includegraphics[width=\textwidth]{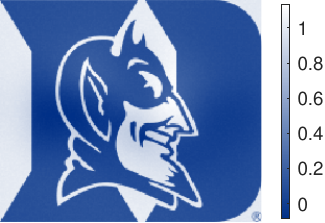}
        \end{subfigure}
        \begin{subfigure}[t]{0.32\textwidth}
            \flushleft
            \caption*{Final $q^{(18)}$ rel. err. $1.68\times 10^{-4}$}
            \includegraphics[width=\textwidth]{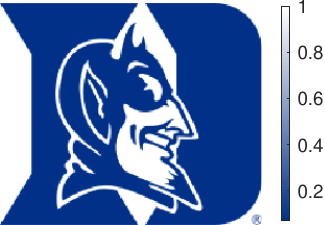}
        \end{subfigure}
        \begin{subfigure}[t]{0.32\textwidth}
            \flushleft
            \caption*{$\phi_0$ rel. err. vs num. outer it.}
            \includegraphics[width=\textwidth]{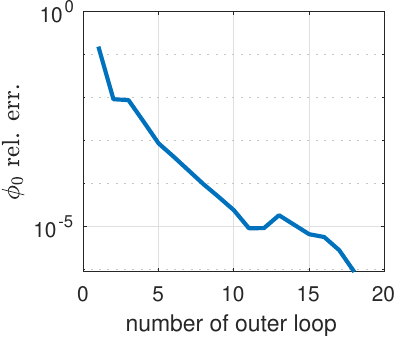}
        \end{subfigure}
    \caption{ECI efficacy with 2D example (\S \ref{subsec: num effective}).
    The measurement-induced term $M(\phi_0):=\frac{1}{T}\int_{\mathbb{T}^d}\phi_0(x)\dd x-\nu\Delta\phi_0+H(\nabla\phi_0)$ retains substantial information about the true potential $q$. 
    The ECI algorithm utilizing $M(\phi_0)$ rapidly recovers the main structure of the true potential. 
    }
    \label{fig: bluedevil}
\end{figure}

The true ambient potential in 1D, shown in Figure \ref{fig: 1dillu} (top left), is defined as:
\begin{equation*}
\begin{aligned}
    q(x) =& (\mathbf{1}_{x<0.4}+\mathbf{1}_{x>0.7})(\sin(20\pi x) \, e^{-10(x-0.5)^2}) 
    + \mathbf{1}_{0.4<x<0.7}(-e^{x}) 
    + \mathbf{1}_{x>0.7}(0.2\sin(100\pi x)) \\
    &- \mathbf{1}_{0.3<x<0.35}
    + \mathbf{1}_{0.6<x<0.65},
\end{aligned}
\end{equation*}
which contains both multiscale and discontinuous features.
The implementation of PIT from \cite{ren2024policy} relies on an optimization solver, and as shown in Figure~\ref{fig: 1dillu} (top left), it struggles to capture the high-frequency oscillations and discontinuities present in the true potential.
The true potential in 2D, shown in Figure \ref{fig: bluedevil} (top left), is set as the ``blue devil'' logo of Duke University.

For the inverse solver, the initial estimate $q^{(0)}(x)$ is chosen randomly from a uniform distribution $U(a,b)$, where $a=\min_x q(x)$ and $b=\max_x q(x)$. 
In each ECI iteration, the forward problem is solved to an accuracy of $10^{-8}$. 
The inverse solver terminates when the measurement relative error falls below $10^{-9}$ for the 1D example and $10^{-6}$ for the 2D example.

Figures \ref{fig: 1dillu} and \ref{fig: bluedevil} display the true potential $q$, the measurement-induced term in ECI
\begin{equation}
\label{eq:measurement}
M(\phi_0):=\frac{1}{T}\int_{\mathbb{T}^d}\phi_0(x)\dd x-\nu\Delta\phi_0+H(\nabla\phi_0),
\end{equation}
the ECI result after one iteration $q^{(1)}$, the final ECI reconstruction $q^{(K)}$, and the convergence plot of ECI.

From the true potential $q$ and measurement-induced term $M(\phi_0)$, we observe that even with a relatively large viscosity, the measurement still contains significant information about the fine detail of the potential. 
ECI utilizes the information of $M(\phi_0)$ and therefore rapidly captures the main structure of the true potential after just one iteration and converges to a highly accurate solution within 30 iterations.
Quantitatively, the 1D example converges in 23 iterations with a measurement relative error of $5.30\times 10^{-10}$ and an ambient potential relative error of $2.88\times 10^{-9}$. 
The 2D example converges in 18 iterations with a measurement relative error of $8.94\times 10^{-7}$ and an ambient potential relative error of $1.68\times 10^{-4}$.

\subsection{Scalability of ECI}
\label{subsec: num comp ECI policy}

\begin{figure}[htb]
    \centering
    \includegraphics[height=4.2cm]{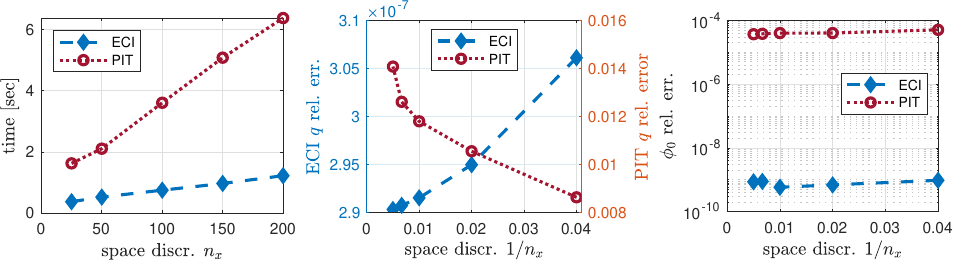}
    \caption{Scalability with respect to grid (\S \ref{subsec: num comp ECI policy}).
    The computational time is approximately linear in $n_x$, and the relative error improves as $n_x$ increases.
    }
    \label{fig: comp eci vs grid}
\end{figure}

This section investigates the scalability of the ECI algorithm and compares its performance with the policy iteration (PIT) method proposed and implemented in \cite{ren2024policy}. The test setting follows the same configuration as Section 5.1 in \cite{ren2024policy}.
We consider the MFG spatially on $[0,1]$ with periodic boundary conditions and $T=1$. 
In the discretization, we fix $\Delta t = 0.01 ~(n_t=100)$ and varies $\Delta x$.
The viscosity is $\nu=0.3$, and the Hamiltonian is $H(p)=\frac{1}{2}|p|^2$. The initial distribution is Gaussian, $\rho_0(x) = \rho_G(x;0.5,\frac{1}{4\sqrt{5}})$, the interaction cost is $f(x,\rho)=\rho(x)^2$, and the terminal cost is independent of the distribution, $f_T(x) = -\rho_0(x)$.
The true ambient potential is set as:
$$q(x)=0.1( \sin{(2\pi x - \sin{(4\pi x)})} + \exp{(\cos{(2\pi x)})} ).$$
For the inverse solver, the initial guess is $q^{(0)}(x)=0.15$. For ECI, in each iteration, we run the forward solver to an accuracy of $10^{-6}$. 
The inverse solver terminates when the measurement relative error is less than $10^{-9}$.

Figure~\ref{fig: comp eci vs grid} presents the computational time versus spatial discretization (left, with $n_x = 25, 50, 100, 150, 200$), and the relative errors of $q^{(K)}$ (center) and $\phi_0^{(K)}$ (right) plotted against the grid step size $\Delta_x = \frac{1}{n_x}$.
Our ECI method demonstrates high computational efficiency, requiring approximately 2 seconds for $n_x = 200$, whereas PIT takes about 6 seconds for the same resolution. Moreover, ECI achieves significantly lower relative errors: for the potential $q$, our method maintains errors around $10^{-7}$ across all grids, compared to approximately $10^{-2}$ by PIT; for the measurement $\phi_0$, ECI achieves errors around $10^{-8}$, while PIT reaches only $10^{-4}$.
In addition, our method benefits from grid refinement: the relative error in recovering the potential $q$ decreases as the grid is refined, which is not observed for the policy iteration approach. This highlights a key advantage of ECI over methods that rely on optimization solvers, where optimization error can dominate and limit the overall accuracy, regardless of the discretization.

While the computational cost and accuracy of our ECI method scale well with respect to grid size, we acknowledge that it is sensitive to noise compared to the optimization-based implementation of policy iteration \cite{ren2024policy}. Developing robust adaptations of ECI to handle noisy measurements remains an important direction for future research.

\subsection{Comparison of ECI and BRI}
\label{subsec: num comp}

\begin{figure}[htb]
    \centering
    \includegraphics[width=\textwidth]{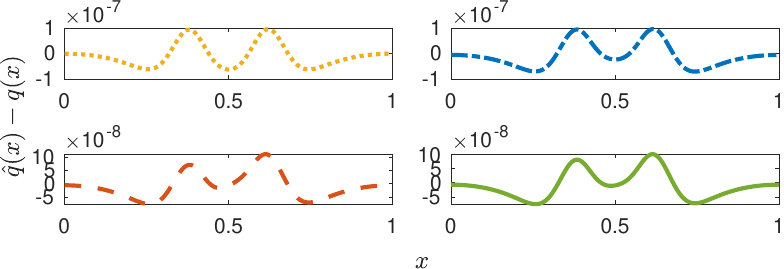}
    \caption{Comparison of ECI and BRI (\S \ref{subsec: num comp}).
    The plots show the errors of recovered potential $q^{(K)} - q$ for BRI1(0.2) (top left), BRI1(0.5) (top right), BRI5(0.5) (bottom left), and ECI (bottom right). 
    Starting from the same initial guess $q^{(0)} = 0$, all methods recover the true potential pointwise with a relative error less than $10^{-6}$.
    }
    \label{fig: comp eci bri q}
\end{figure}

\begin{figure}[htb]
    \flushleft
    \includegraphics[width=0.96\textwidth]{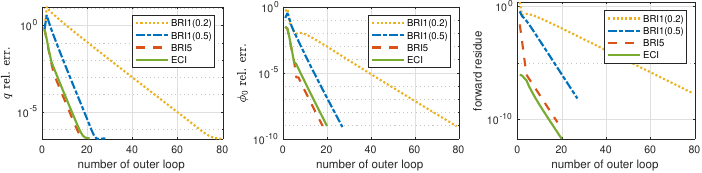}
    \includegraphics[width=0.64\textwidth]{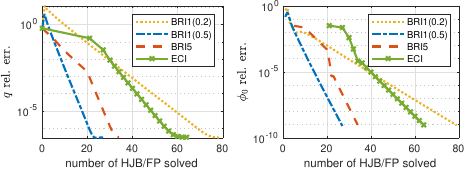}
    \caption{Comparison of ECI and BRI (\S \ref{subsec: num comp}).
    The top row shows the relative error of the potential $q^{(k)}$ (left) and of the measurement $\phi_0^{(k)}$ (center) and the residue of the forward solver (right) versus the number of outer loops.
    The bottom row incorporates the cost of the inner loop iterations and plots the relative error of the potential $q^{(k)}$ (left) and the measurement $\phi_0^{(k)}$ (center) versus the number of HJB/FP solves. Each marker of ECI corresponds to an outer loop iteration.
    BRI1 and BRI5 achieve a similar relative error as ECI in a comparable number of outer loops but with fewer calls to the HJB and FP solvers.
    }
    \label{fig: comp eci bri err}
\end{figure}

This section compares the performance of ECI and BRI and examines how the choice of the weight $\delta_{k,n}$ and the number of inner-loop $n$ affect the convergence of BRI.
Numerical results demonstrate that in these double-loop algorithms, the accuracy of both the inner-loop forward solution and outer-loop solution improves with the outer-loop iterations, even when the inner-loop best-response iteration is limited to a fixed number. Consequently, BRI is more computationally efficient than ECI.

The test setting is the same as Section \ref{subsec: num comp ECI policy}, except that this section uses a finer discretization with $\Delta x = 10^{-3}$ and $\Delta t = 2 \times 10^{-3}$.
For the inverse solver, the initial guess and the ECI settings remain the same as in Section \ref{subsec: num comp ECI policy}. For BRI, we run the forward solver with $N = 1$, $\delta_{k,n} = 0.2$; $N = 1$, $\delta_{k,n} = 0.5$; and $N = 5$, $\delta_{k,n} = 0.5$, labeling the outputs as BRI1(0.2), BRI1(0.5), and BRI5, respectively.
The inverse solver terminates when the relative error of the measurement falls below $10^{-9}$.

As shown in Figures \ref{fig: comp eci bri q}, ECI, BRI1(0.2), BRI1(0.5), and BRI5 all recover the true ambient potential with a pointwise absolute error of less than $10^{-6}$.
In Figure \ref{fig: comp eci bri err}, the top row plot the relative error of the potential $q^{(k)}$ (left) and the measurement $\phi_0^{(k)}$ (center) versus the number of outer loops, as well as forward solver residue (right), i.e., the residue of $(\rho^{(k)},\phi^{(k)})$ to the MFG system with $q^{(k)}$.
The right panel shows that even for BRI1 and BRI5, where the inner loop only runs for a fixed number of iterations, the forward solver residue decreases as the outer loop progresses. As a result, the convergence rates of BRI5 and ECI with respect to the number of outer loops are nearly identical.
The right panel also shows that when $\delta_{k,n}$ is small, the forward residue decays more slowly. As a result, BRI1(0.2) converges more slowly than BRI1(0.5) in the left and center panels. 
The bottom row of Figure \ref{fig: comp eci bri err} incorporates the cost of the inner loop iterations and plots the relative error of the potential $q^{(k)}$ (left) and the measurement $\phi_0^{(k)}$ (center) versus the number of HJB/FP solves.
Each marker for ECI corresponds to an outer loop. Initially, it takes several iterations for the forward solver to converge, but later, the forward solver converges in a single iteration.
However, when the guess of the potential $q^{(k)}$ is not very accurate at the beginning, more inner loop iterations do not significantly improve the accuracy of the potential. As a result, BRI1 and BRI5 achieve the same accuracy as ECI with fewer calls to the HJB and FP solvers, making them more efficient.

\subsection{Remarks on BRI}
\label{subsec: num bri}

\begin{figure}[htb]
    \flushleft
    \includegraphics[width=0.64\textwidth]{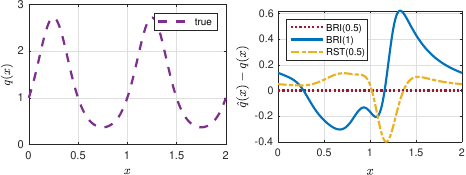}
    \includegraphics[width=0.96\textwidth]{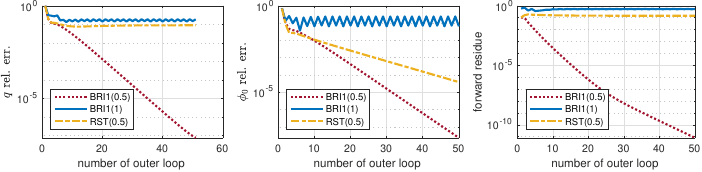}
    \caption{Remarks on BRI (\S \ref{subsec: num bri}).
    The top row shows the true potential $q$ (left) and the errors $q^{(K)}-q$ for BRI(0.5), BRI(1) and RST(0.5) (right).
    The bottom row shows the relative error of the potential $q^{(k)}$ (left) and of the measurement $\phi_0^{(k)}$ (center) and the residue of the forward solver (right) versus the number of outer loops.
    BRI(1) and RST(0.5) are unable to improve forward solver accuracy and therefore fail to recover the true obstacle. 
    }
    \label{fig: bri rmk}
\end{figure}

This section demonstrates that it is necessary to initialize the inner loop of BRI with the previous distribution $\tilde{\rho}^{(k,0)}=\rho^{(k-1)}$ and to choose the weight $\delta_{k,n}$ properly.

We choose the space domain to be $[0,2]$ with periodic boundary condition, $T=1,\nu=0.1$ and the Hamiltonian to be $H(p)=\frac{1}{2}|p|^2$. The initial distribution is Gaussian with mean $1$ and variance $0.04$. The interaction cost is $f(x,\rho)=\rho(x)$ and the terminal cost is absent $f_T(x,\rho)=0$.
The true potential $q(x)=\exp{(\sin{(2\pi x)})}$.

For all tests, we initialize the algorithms with $q^{(0)} = 0$. We run BRI (Alg.~\ref{alg: BRI}) with $N = 1$ and $\delta_{k,n} = 0.5$, and label the result as BRI(0.5). We also run BRI with $N = 1$ and $\delta_{k,n} = 1$, labeled as BRI(1).
As discussed in Section~\ref{subsec: bri}, the initialization of the inner loop is also important. In BRI, we set $\tilde{\rho}^{(k,0)} = \rho^{(k-1)}$. In this numerical experiment, we additionally test using $\tilde{\rho}^{(k,0)} = \rho^{(0)}$, where $\rho^{(0)}_t = \rho_0$ is a static flow. All other settings remain the same as in BRI, and we choose $\delta_{k,n} = 0.5$. We refer to this variant as RST(0.5).

As shown in Figure~\ref{fig: bri rmk} (top center and bottom left), both BRI(1) and RST(0.5) fail to recover the true potential.
For BRI(1), the choice of $\delta_{k,n} = 1$ is too aggressive. As a result, the forward residual does not decay (bottom right), and both the estimated potential and the measurement exhibit oscillatory behavior (bottom left and bottom center). Similar oscillations have been observed in the forward solver fictitious play when using $\delta_{k,n} = 1$, as reported in \cite{yu2024ficplay}.
In the case of RST(0.5), although $\delta_{k,n} = 0.5$ is a moderate choice, the failure arises from the initialization of the inner loop. Instead of using the previous approximation $\rho^{(k-1)}$, RST(0.5) initializes with a fixed static flow $\rho^{(0)}$, independent of the current iteration. As a result, a fixed number of fictitious play iterations only approximates an intermediate stage of the learning dynamics, rather than the Nash equilibrium. Consequently, the outer-loop update matches the intermediate value function with the observed measurement, which leads to a decrease in the measurement error (bottom center). However, since the approximation does not correspond to a Nash equilibrium, the method fails to recover the true potential (bottom left).

\subsection{Acceleration with Hierarchical Grid}
\label{subsec: num acc}

This section demonstrates that the hierarchical grid approach can accelerate the detection of the low-frequency component of the ambient potential.

We consider the two different potentials $q$, one with a low-frequency component and the other with a high-frequency component:
\begin{eqnarray*}
    \text{Panel 1: } & q(x) = e^x \sin{(2\pi x)}.\\
    \text{Panel 2: } & q(x) = e^x \sin{(8\pi x)}.
\end{eqnarray*}

We choose the space domain to be $[-1,1]$ with periodic boundary condition, $T=1,\nu=0.1$, and the Hamiltonian to be $H(p)=\frac{1}{2}|p|^2$. The initial distribution is uniform, the interaction cost is $f(x,\rho)=\rho(x)$ and the terminal cost is absent $f_T(x,\rho)=0$.

\begin{figure}[htb]
    \centering
    \includegraphics[width=\textwidth]{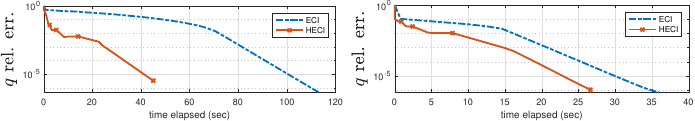}
    \caption{Acceleration with hierarchical grid (\S \ref{subsec: num acc}).
    Relative error of $q$ versus the time elapsed of ECI and HECI for $q(x)=e^x\sin{(2\pi x)}$ (left) and $q(x)=e^x\sin{(8\pi x)}$ (right). 
    The markers for HECI correspond to the convergence on different levels of the hierarchical grid.
    HECI greatly accelerates the convergence for low-frequency potentials.
    }
    \label{fig: comp eci heci}
\end{figure}

We initialize the algorithm with $q^{(0)}(x)=0$ and run the ECI algorithm on a uniform grid with $n_x=2^{10}$ and $n_t=2^{9}$.
We also run the hierarchical ECI (HECI) algorithm on a hierarchical grid with $L=4$ levels, where the coarse grid $\calG_0$ has $n_{x,0}=2^7$ and $n_{t,0}=2^6$, and each subsequent grid $\calG_l$ is obtained by refining the previous grid by a factor of 2 in both space and time.
Figure \ref{fig: comp eci heci} shows the relative error of the estimated ambient potential versus the time elapsed for ECI and HECI for the tests. The markers for HECI correspond to the convergence on different levels of the hierarchical grid.

Quantitatively, for the low-frequency potential, ECI requires 77 outer iterations and 276 HJB/FP solves to converge, whereas for the high-frequency potential, it converges in just 32 outer iterations and 91 HJB/FP solves. 
We observe that for the low-frequency potential, HECI achieves a $10^{-1}$ accuracy on the coarse grid in less than 10 seconds, while ECI requires more than 60 seconds to reach the same accuracy on the fine grid. As a result, HECI greatly accelerates the convergence for low-frequency potentials.
And for the high-frequency potential, ECI takes very few iterations to reach a high accuracy on the fine grid, and HECI does not provide a significant acceleration.

In summary, ECI rapidly captures the high-frequency components of the potential, while HECI provides significant acceleration in recovering the low-frequency components.

\subsection{Numerical Evidences of Convergence Mechanism}
\label{subsec: num cvg}

In the end, we conduct numerical experiments to gain a preliminary understanding of the inverse problem and the convergence of the algorithm.

We choose the spatial domain to be $[-1,1]$ with periodic boundary conditions, and the Hamiltonian is $H(p)=\frac{1}{2}|p|^2$. The initial distribution is $\rho_0(x)=\rho_G(x;0,0.2)$, and the terminal cost is absent, i.e., $f_T(x,\rho)=0$.
The ground truth ambient potential is set as 
$$q(x)=0.1(\exp{(\sin{(2\pi x^3)})} + (x+1)(x-1)(x-0.5)-2),$$ 
and the algorithm is initialized with $q^{(0)}(x)=0$.

In the forward problem, both the viscosity and the monotonicity of the interaction and terminal costs play significant roles in well-posedness and algorithmic convergence. For a Hamiltonian that is strictly convex in $p$, monotone interaction and terminal costs guarantee uniqueness of the solution to the forward system. Moreover, fictitious play is provably convergent for monotone costs and is observed to converge more rapidly as viscosity increases. 
Additionally, due to the connection between the MFG system and semilinear parabolic equations, the terminal time $T$ may also influence the problem. In parabolic equations, both viscosity and time contribute to diffusion, so the product $\nu T$ is a key parameter determining the state at time $T$.

To numerically investigate the inverse problem and the algorithm, we vary the viscosity $\nu$ and terminal time $T$, as detailed in Table \ref{tab: test cvg mono} for a monotone interaction cost $f(x,\rho)=\rho(x)$ and in Table \ref{tab: test cvg nonmono} for a non-monotone interaction cost $f(x,\rho)=-\rho(x)$. 

\begin{figure}[htb]
    \centering
    \includegraphics[width=\textwidth]{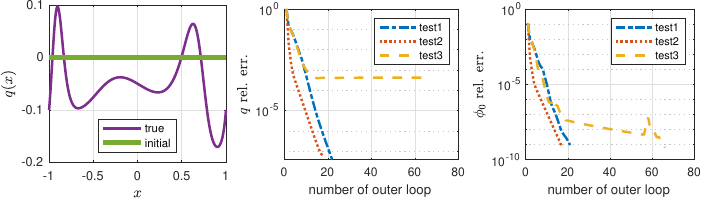}
    \includegraphics[width=\textwidth]{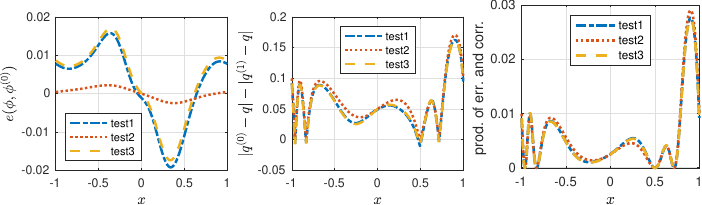}
    \caption{Convergence study with monotone interaction cost $f(x,\rho)=\rho(x)$ (\S \ref{subsec: num cvg}, Table \ref{tab: test cvg mono}). Test 1 is the base case with $\nu=0.1$ and $T=1$. Test 2 with $\nu=0.5$ and $T=1$ explores the effect of viscosity, and Test 3 with $\nu=0.1$ and $T=5$ explores the effect of terminal time.
    }
    \label{fig: test cvg mono}
\end{figure}

\begin{figure}[htb]
    \centering
    \includegraphics[width=\textwidth]{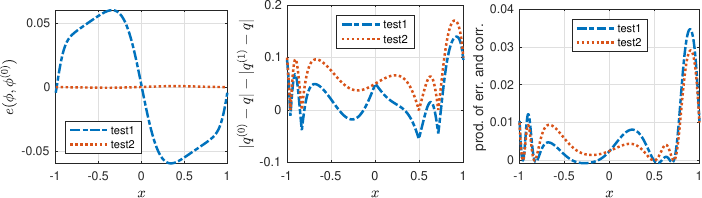}
    \caption{Convergence study with non-monotone interaction cost $f(x,\rho)=-\rho(x)$ 
    (\S \ref{subsec: num cvg}, Table \ref{tab: test cvg nonmono}). Test 4  
    with $\nu=0.1$ diverges and test 5  
    with $\nu=1$ converges in 13 iterations with a measurement relative error of $8.86\times 10^{-10}$.
    }
    \label{fig: test cvg nonmono}
\end{figure}

For each setting, we conduct two types of tests. First, we run the inverse solver using the measurement $\phi_0$ and initialization $q^{(0)}$ to compare convergence rates under different parameter choices and report them in Tables \ref{tab: test cvg mono} and \ref{tab: test cvg nonmono}. Second, we run the forward solver with both the true potential $q$ and the initial guess $q^{(0)}$ to obtain the corresponding measurements $\phi_0$ and $\phi_0^{(0)}$. 
Recall equations \eqref{eq: update correction term}, \eqref{eq: update in terms of correction}, and \eqref{eq: update error term}: this setup allows us to compute the correction $c(q,q^{(0)})$, the error in correction $e(\phi,\phi^{(0)})$, and the product of error and correction (PEC), which are important indicators of algorithm performance.
More precisely, the error in correction $e(\phi,\phi^{(k)})$ defined in \eqref{eq: update error term} is also the error of the updated estimate $q^{(k+1)}$. Thus, the ECI update exactly recovers the ground truth when $e(\phi,\phi^{(k)})$ is zero pointwise. If $|e(\phi,\phi^{(0)})| < |q^{(0)} - q|$ pointwise, the error decreases at each state. The product of error and correction, $\mathrm{PEC} = (q - q^{(0)})c(q,q^{(0)})$, is positive when the correction is in the correct direction. The plots are in Figures \ref{fig: test cvg mono} and \ref{fig: test cvg nonmono}.

\begin{table}[htb]
    \centering
    \caption{Tests with monotone interaction cost $f(x,\rho)=\rho(x)$.}
    \begin{tabular}{c|ccccc}
        \toprule
         & $\nu$ & $T$ & cvg. it. & rel. err. of $\phi_0$ & rel. err. of $q$ \\
        \midrule
        Test 1 & 0.1 & 1  & 21 & $9.35\times 10^{-10}$ & $3.90\times 10^{-8}$ \\ \hline
        Test 2 & 0.5 & 1  & 17 & $8.20\times 10^{-10}$ & $5.80\times 10^{-8}$ \\ \hline
        Test 3 & 0.1 & 5  & 66 & $6.30\times 10^{-10}$ & $4.22\times 10^{-4}$ \\ 
        \bottomrule
    \end{tabular}
    \label{tab: test cvg mono}
\end{table}
\begin{table}[htb]
    \centering
    \caption{Tests with nonmonotone interaction cost $f(x,\rho)=-\rho(x)$.}
    \begin{tabular}{c|ccccc}
        \toprule
         & $\nu$ & $T$ & cvg. it. & rel. err of $\phi_0$ & rel. err. of $q$ \\
        \midrule
        Test 4 & 0.1 & 1  & N/A & 1.09 & 19.02 \\ \hline
        Test 5 & 1   & 1  & 13  & $8.86\times10^{-10}$ & $2.41\times 10^{-6}$ \\ \hline
        \bottomrule
    \end{tabular}
    \label{tab: test cvg nonmono}
\end{table}

\begin{figure}[htb]
    \flushleft
    \includegraphics[width=0.96\textwidth]{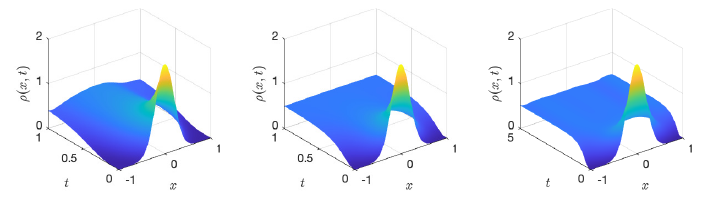}
    \includegraphics[width=0.64\textwidth]{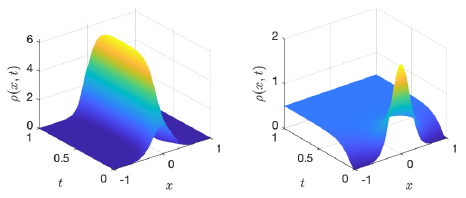}
    \caption{Convergence study (\S \ref{subsec: num cvg}). 
    Ground truth $\rho$ for Table \ref{tab: test cvg mono}, Test 1-3 (top row, left to right) and Table \ref{tab: test cvg nonmono}, Test 4-5 (bottom row, left to right).
    For Test 4 (bottom left, non-monotone cost), the population density remains highly localized. This limits the information content of the measurement, making the inverse problem intrinsically challenging in this regime.
    }  
    \label{fig: test cvg rhotru}
\end{figure}

For the monotone interaction cost, the algorithm is effective in all three tests and suggests that large viscosity $\nu$ is beneficial for convergence while larger terminal time $T$ may lead to ill-posedness of the inverse problem.
Figure \ref{fig: test cvg mono} (bottom left) shows that the error in correction is small across all three tests, particularly when the viscosity $\nu$ is large. This indicates that the first ECI update is highly effective. The bottom center and right panels further demonstrate that, for most states $x$, the ECI update moves in the correct direction and yields pointwise improvement, though some overshooting may occur. The effectiveness of the algorithm is also evident in the convergence results in Table \ref{tab: test cvg mono}, where the measurement relative error for $\phi_0$ consistently reaches around $10^{-9}$. However, both Table \ref{tab: test cvg mono} and Figure \ref{fig: test cvg rhotru} (top center and right panels) reveal that for Test 3 ($T=5$), the relative error in the recovered potential plateaus at a larger value compared to Tests 1 and 2. This suggests that the inverse problem becomes ill-posed when the terminal time $T$ is large, since different potentials may produce nearly indistinguishable measurements $\phi_0$.

The inverse problem with a non-monotone interaction cost is more challenging; however, increasing the viscosity $\nu$ can improve stability and convergence. Table \ref{tab: test cvg nonmono} shows that the algorithm fails to converge in Test 4 ($\nu=0.1$), while for larger viscosity (Test 5, $\nu=1$), it converges rapidly, achieving a measurement error of $10^{-9}$. However, the potential error ($10^{-5}$) is still higher than what is observed for monotone costs ($10^{-7}$). Figure \ref{fig: test cvg nonmono} (right panel) indicates that the correction is generally in the correct direction for most states, but the algorithm tends to overshoot in many locations when the viscosity is small (center panel), leading to divergence in Test 4.

Heuristically, the failure of the algorithm in Test 4 (non-monotone cost) appears to stem from the intrinsic ill-posedness of the inverse problem in this regime. Figure \ref{fig: test cvg rhotru} shows the ground truth $\rho$ for the tests in Tables \ref{tab: test cvg mono} and \ref{tab: test cvg nonmono}. In Test 4 (bottom left, non-monotone cost), the population density $\rho$ remains highly localized, with a peak near $x=0$ and negligible values away. This contrasts with the other tests, where the density is more broadly distributed across the spatial domain. This implies that the population does not sufficiently explore the entire domain, limiting the information content of the measurement $\phi_0$ and making the inverse problem intrinsically difficult in this setting.

\section{Conclusion and Future Work}
\label{sec: conclusion}

In this paper, we propose simple and constructive algorithms for recovering the ambient potential in mean-field games (MFGs) from measurements of the value function at the Nash equilibrium and demonstrate their effectiveness through a range of numerical examples. 
While our results highlight the practical potential of these methods, establishing the well-posedness of the inverse problem and proving the convergence of the algorithms remain important directions for future research. Numerical experiments and the connection between our inverse MFG problem and classical inverse problems for linear parabolic equations provide valuable insights and suggest promising avenues for further theoretical investigation. 
However, how to adapt ECI to recover the potential from noisy observations remains an open problem.

\section*{Acknowledgments}
We sincerely thank Prof. Kui Ren and Prof. Shanyin Tong for generously sharing the code associated with \cite{ren2024policy}. 


\bibliography{reference}
\bibliographystyle{siam}

\end{document}